\documentclass[a4paper,12pt]{article}

\newtheorem{remark}{Remark}
\newtheorem{problem}{Problem}

\usepackage{amsmath}
\usepackage{amssymb}
\usepackage{amsfonts}
\usepackage{graphicx}
\usepackage{mathtools}

\usepackage{xspace}

 \let\forallOLD\forall \renewcommand{\forall}{\ensuremath{\forallOLD\,}}
\newcommand{\iset}[2]{#1,\ldots,#2}

\newcommand{\intT}{\ensuremath{\int_{\mathcal{T}}}}

\usepackage{enumerate}
\newcounter{myAlgoCounter} 

\title{\bf 
  Estimation of consistent parameter sets for continuous-time nonlinear
  systems using occupation measures and LMI relaxations}

\begin{document}

\author{Stefan Streif$^1$, Philipp Rumschinski$^1$,\\ Didier Henrion$^{2,3,4}$ and Rolf Findeisen$^1$}

\footnotetext[1]{Laboratory for Systems Theory and Automatic Control, Institute for Automation Engineering,
Otto-von-Guericke-University~Magdeburg, Germany.}
\footnotetext[2]{CNRS, LAAS, University of Toulouse, France.}
\footnotetext[3]{Faculty of Electrical Engineering, Czech Technical University in Prague,
Czech Republic}

\date{Draft of \today}

\maketitle

\begin{abstract}
    Obtaining initial conditions and parameterizations leading to a model consistent with available measurements or
    safety specifications is important for many applications. Examples include model (in-)validation, prediction, fault
    diagnosis, and controller design.  We present an approach to determine inner- and outer-approximations of the set
    containing all consistent initial conditions/parameterizations for nonlinear continuous-time systems.  These
    approximations are found by occupation measures that encode the system dynamics and measurements, and give rise to
    an infinite-dimensional linear program.  We exploit the flexibility and linearity of the decision problem to
    incorporate uncertain-but-bounded and pointwise-in-time state and output constraints, a feature which was not
    addressed in previous works.  The infinite-dimensional linear program is relaxed by a hierarchy of LMI problems that
    provide certificates in case no consistent initial condition/parameterization exists.  Furthermore, the applied LMI
    relaxation guarantees that the approximations converge (almost uniformly) to the true consistent set.  We illustrate
    the approach with a biochemical reaction network involving unknown initial conditions and parameters.
\end{abstract}

\section{Introduction}

The computation of guaranteed inner- and outer-approximations of consistent parameter sets of uncertain dynamical
systems is important for many applications including model-based analysis and verification, system identification, model
(in-)validation and controller design
\cite{Milanese_1991_Automatica__Estimation_theory_set_membership,Anderson_Papachristodoulou_2009_BMC__On_validation_invalidation_biological_models,Prajna_Rantzer_2007_SIAM__Convex_programs_temporal_verifications_nonlinear_systems,Prajna_2006_Automatica__Barrier_certificate_model_validation}.

We consider the derivation of such approximations for polynomial systems subject to unknown-but-bounded (or
error-in-variables) data, e.\,g. measurements.  Different methods are available in literature to address this problem.
For discrete-time systems, for instance interval analysis methods, e.\,g.
\cite{Jaulin_etAl_2001_Interval_analysis_parameter_state_estimation}, or relaxation-based methods, e.\,g.
\cite{Streif_etAl_2012_ADMIT,Rumschinski_etAl_2010_Set-based_parameter_estimation_model_invalidation,Borchers_etAl_2009_CDC__Set-based_estimation,Cerone_etAl_2012__Set-membership_identification_convex_relaxation}
can be employed.  However, both approaches are not directly applicable to continuous-time systems without additional
assumptions. For instance in \cite{Fey_Bullinger_2010_SDP} it was assumed that the time-derivatives of the states are
available as measurements, therefore, resulting in a steady-state problem similar to \cite{Kupfer:2007vl}.

One possibility to address continuous-time systems with the mentioned methods is by discrete-time approximations,
e.\,g. obtained by numerical integration.  However, due to the
discretization error the consistent parameter sets of
continuous-time and discrete-time model do not necessarily overlap and, thus, wrong conclusions with respect to model
validity are possible \cite{Rumschinski:2010tx}.  A common approach to limit the discretization error relies on
higher-order Taylor approximations resulting in verified integration methods
\cite{Raissi_etAl_2004_Set-membership_state_parameter_estimation,Johnson_Tucker_2008_Parameter_reconstruction_differential_equation_noisy_data,Nedialkov_etAl_1999_Validated_solutions_initial_value_problems_ODEs,Lin_Stadtherr_2007_Validated_solutions_IVP_ODEs,Berz_Makino_1998_Verified_integration_ODEs,Rauh_etAl_2009_Interval_arithmetic_solving_DAEs},
but again the results typically depend on the discretization error.

A more direct approach uses McCormick relaxations or differential
inequalities, see
\cite{Scott_Barton_2013_Relaxation_parametric_solutions_ODEs_differential_inequalities}
and references therein). However,
deriving the needed tight state bounds can be difficult.

Methods using barrier certificates and sum-of-squares (SOS) polynomials
\cite{Prajna_Rantzer_2007_SIAM__Convex_programs_temporal_verifications_nonlinear_systems,Prajna_2006_Automatica__Barrier_certificate_model_validation,Anderson_Papachristodoulou_2009_BMC__On_validation_invalidation_biological_models}
allow the continuous-time dynamics to be considered directly without
numerical integration. However, only few converse results for these,
i.\,e.\ existence of barrier certificates, are known.  Furthermore, to the best of our knowledge, no results with
respect to approximations of consistent parameter sets exist so far.

The main contribution of this work is the use of occupation measures \cite{Mezic_Runolfsson_2008_Uncertainty_propagation} to derive
guaranteed polynomial inner- and outer-approximations of the consistent parameter sets for continuous-time nonlinear systems based
on results presented in
\cite{Henrion_Korda_2013_ROA_occupation_measure,Korda_etAl_2013__Inner_approximations_occupation_measure,Lasserre_2009}.
The reformulation in terms of occupation measures leads to a linear but infinite-dimensional decision problem. However,
its relaxation using truncated moment matrices and their dual SOS
polynomials is a finite-dimensional convex
problem.  One particular feature of the employed relaxation is the (almost uniform) asymptotic convergence of the
approximations to the true consistent parameter set.  Another advantage is that the continuous-time dynamics are
completely encoded in the decision problem and, thus, no numerical integration is necessary.  Furthermore, we exploit the
flexibility and linearity of the decision problem to incorporate uncertain-but-bounded and pointwise-in-time state and
output constraints, a feature which was not addressed in previous works.

This contribution is structured as follows. In Section~\ref{sec:problem} we formalize the problem setup, in particular,
the considered system class, the description of the uncertain data and the desired properties of the approximations.  To
obtain constraints for a convex optimization problem, we reformulate in Section~\ref{sec:occup_meas} the polynomial
continuous-time dynamics in terms of occupation measures.
In
Section~\ref{sec:appli}, we adapt this approach to the consistent parameter
estimation problem, which is reformulated as an
infinite-dimensional linear programming problem. Its solution is approached numerically with a hierarchy of
finite-dimensional semi-definite programs.  We show how to derive outer- and inner-approximations, as well as
certificates of inconsistency.  The approach is demonstrated in Section~\ref{sec:example} while its advantages and
computational issues are discussed in Section~\ref{sec:conclusion}.

\section{Problem Formulation}
\label{sec:problem}

In this section we state the problem of set-based parameter estimation for continuous-time nonlinear systems of the
following form
\begin{subequations}
    \label{eq:system}
    \begin{align}
        \dot{x}(t) &= f\bigl(t,x(t)\bigr), \qquad x(t_0) = x_0, \label{eq:system:ode}\\
        y(t) &= h\bigl(x(t)\bigr).  \label{eq:system:output}
    \end{align}
\end{subequations}
Here $t \in [0,1]$ is the time, the states are denoted by $x \in \mathbb{R}^{n_x}$ and the outputs by $y \in
\mathbb{R}^{n_y}$.  Note that the terminal time is set to one without loss of generality, after a suitable time-scaling
of the dynamics.  Initial conditions (at $t=0$) are denoted by $x_0$.  We assume the vector field $f:
\mathbb{R}\times\mathbb{R}^{n_x}\rightarrow \mathbb{R}^{n_x}$ and $h: \mathbb{R}\times\mathbb{R}^{n_x}\rightarrow
\mathbb{R}^{n_y}$ to be polynomial maps.  To simplify notation, we denote the state vector at time $t_k\in[0,1]$ by $x_k
= x(t_k)$.

Note that time-invariant parameters $p \in \mathbb{R}^{n_p}$ can be included in (\ref{eq:system}) by defining states
with constant dynamics $\forall i \in \iset{1}{n_p}$:
\begin{equation}
    \label{eq:system_p}
        \dot{x}_{i} = f_{i}\bigl(t,x(t)\bigr) = 0,\quad
        \text{where }
        p_i := x_{i}(0).
\end{equation}

That is, the state vector $x$ contains $n_{x}-n_p$ true dynamic variables and $n_p$ constant parameters.  The parameter
values are then given by the initial conditions of the corresponding variables $x_i$. This formulation unifies the tasks
of initial condition and parameter estimation, and also simplifies the notation and analysis using occupation measures
and moment relaxations (see next sections).

We assume that constraints on the states, output measurements and initial conditions (including the parameters) are
given by polynomial inequalities $g(\cdot) \geq 0$. We rewrite the output constraints $g_y\bigl(y(t)\bigr) \geq 0$ as state
constraints using the polynomial output map, i.\,e.  $g_{y}\bigl(h(t, x)\bigr) \geq 0$.

Let $m_x$ time-independent constraints on the states be given in the form $x(t) \in \mathcal{X}, \forall t \in
\left[0,1\right]$ with
\begin{align}
    \label{eq:data_x}
    \mathcal{X} \coloneqq \Bigl\{ x :\ g_{x,i}\bigl(x\bigr) \geq 0, \ \forall i = \iset{1}{m_x} \Bigr\} \subset
    \mathbb{R}^{n_x}.
\end{align}

We additionally assume that we have a finite set of $m_t$ distinct measurements at time points $t_k, k =
\iset{0}{m_t-1}$, such that: $0 = t_0 < t_1 < \cdots < t_{m_t-1} = 1$.
These time points include both the measurements of the output $y(t_k)$ and a priori information on the set of parameters
and initial conditions.  For each time point $t_k$, let the constraints on the states be given in the form $x_k \in
\mathcal{X}_k, k=\iset{0}{m_t-1}$ with
\begin{align}
    \label{eq:data_x_k}
    \mathcal{X}_k \coloneqq \Bigl\{ x :\ g_{x_k,i}\bigl(x\bigr) \geq 0,\ \forall i=\iset{1}{m_{x_k}} \Bigr\} \subseteq
    \mathcal{X}
\end{align}
Using this information, we want to estimate the consistent initial conditions and parameters.

In the following, we denote the admissible state trajectory (an absolutely continuous function of time) of
\eqref{eq:system} starting at fixed $x_0$ with $x(t|x_0)$, i.\,e.
\begin{equation}
    \label{eq:traj_state}
    x(t|x_0) = x_0 + \int_{0}^t f\bigl(\tau,x(\tau)\bigr)\,d\tau.
\end{equation}

Using \eqref{eq:traj_state}, we define the set of consistent initial conditions $\mathcal{X}_0^\ast$ based on the set of consistent
trajectories as follows:
\begin{align}
    \label{eq:consistent:X_0}
    \mathcal{X}_0^\ast \coloneqq \Bigl\{ x_0 : \exists &x(t|x_0) \text{\ s.\,t.\ } \nonumber\\
    &x(t|x_0) \in \mathcal{X},\ \forall t\in [0, T] \text{\ and\ } \\ & x(t_k|x_0) \in \mathcal{X}_k,\ \forall k =
    \iset{0}{m_t-1} \Bigr\}.  \nonumber
\end{align}

With these notations, we can state the following problem:
\begin{problem}[Consistent parameter estimation]\ \\
    \label{prob:param_estim}%
    Find inner-approximations $\mathcal{I}$ 
    and outer-approximations $\mathcal{O}$ 
    of consistent initial conditions (parameters) such that ${\mathcal{I}} \subseteq \mathcal{X}_0^\ast \subseteq
    \mathcal O$.
\end{problem}

For many applications such as model validation and fault detection, it is sufficient to determine if the set
$\mathcal{X}_0^\ast$ is empty or not, we therefore state the following related problem:
\begin{problem}[Certificate of inconsistency]\ \\
    \label{prob:inconsistency}%
    If $\mathcal{X}_0^\ast$ is empty, find a certificate of emptiness.
\end{problem}

Consistency of a single initial condition $x_0$ can easily be checked numerically by solving equation
\eqref{eq:traj_state}. However, to determine the complete set $\mathcal{X}_0^\ast$ is in general difficult, e.\,g. due
to the involved nonlinearities and nonconvexities. We use occupation measures (see Section~\ref{sec:occup_meas}) and
convex relaxations to derive the outer-approximation $\mathcal{O}$ based on
\cite{Henrion_Korda_2013_ROA_occupation_measure}.  The basic idea to find the inner-approximation $\mathcal{I}$ is to
determine guaranteed enclosures of initial conditions that violate some constraints, cf.\ Section~\ref{sec:appli} and
\cite{Korda_etAl_2013__Inner_approximations_occupation_measure}.
Note that the inner-approximation $\mathcal I$ is given by the complement in $\mathcal X$ of the union of these
enclosures.

In the next section, we describe how to deal with the continuous-time dynamics and state constraints using occupation
measures. This enables us to derive convex optimization problems that are used to determine the inner and outer-approximations.

\section{Occupation Measures and Liouville's Equation}
\label{sec:occup_meas}

The crucial idea we employ in this work is to reformulate the parameter estimation problem in terms of occupation
measures.  This has two main advantages. First, it allows us to consider an entire distribution (or measure)
of initial conditions and parameter values and not just single points.  Second, linear
relationships encoding the nonlinear dynamics (i.\,e.\ trajectories) of the system link these initial measures with
corresponding measures at the intermediate time-points for which measurement data are available.


In \cite{Henrion_Korda_2013_ROA_occupation_measure} and the references therein, occupation measures were used to estimate
the region of attraction. Here we first review this approach, before we extend it for parameter estimation in Section
\ref{sec:appli}.  We show that unknown-but-bounded state constraints at intermediate time points, not considered in
\cite{Henrion_Korda_2013_ROA_occupation_measure}, can be handled easily. This is achieved by partitioning the occupation
measure w.r.t. a partitioning of the time interval, cf. Section~\ref{sec:appli}.
To determine the unknown occupation measures, a convex problem is derived. Albeit infinite-dimensional, the convex
problem can be solved efficiently by a hierarchy of finite-dimensional relaxations
\cite{Henrion_Korda_2013_ROA_occupation_measure,Korda_etAl_2013__Inner_approximations_occupation_measure,Lasserre_2009}.

\subsection{Occupation measures}

%
Let $M(\mathcal{A})$ denote the set of finite Borel measures supported on the set $\mathcal A$, which can be interpreted as
elements of the dual space $C(\mathcal{A})'$, i.\,e.\ as bounded linear functionals acting on the set of continuous
functions $C(\mathcal{A})$.  Let ${P}(\mathcal{A})$ denote the set of probability measures on $\mathcal A$, i.\,e.\
those measures $\mu$ of ${M}(\mathcal{A})$ which are nonnegative and normalized to $\mu(\mathcal{A})=1$.

Now assume that the initial condition $x_0$ is not precisely known, but that it can be interpreted as a random variable
whose distribution is described by a probability measure $\mu_0 \in {P}(\mathcal{X})$. We define the occupation measure%
\[ \mu(\mathcal{A}\times\mathcal{B}) := \intT \int_{\mathcal{X}}
I_{\mathcal{A}\times\mathcal{B}}(t,x(t|x_0))\mu_0(dx_0)\,dt \] %
for all subsets $\mathcal{A}\times\mathcal{B}$ in the Borel $\sigma$-algebra of subsets of
$\mathcal{T}\times\mathcal{X}$, where $\mathcal T \subset \mathbb R$ is a time interval and $x(t|x_0)$ is as in
\eqref{eq:traj_state}. Here, $I_{\mathcal{A}}(x)$ is the indicator function of the set $\mathcal A$, which is equal to one if $x \in
\mathcal{A}$, and zero otherwise.

Note that $\mu \in {P}(\mathcal{T}\times\mathcal{X})$ is a probability measure, and that the terminology occupation
measure is motivated by the observation that the value $\intT \mu(dt,\mathcal{B}) = \mu(\mathcal{T}\times\mathcal{B})$
is equal to the total time the trajectory spends in the set $\mathcal{B} \subset \mathcal{X}$.  In addition, note that
$\mu$ encodes the system trajectories, in the sense that if $v \in {C}^{\infty}(\mathcal{T}\times \mathcal{X}; {\mathbb
   R})$ is a smooth test function, and $\mu_0 = \delta_{x_0}$ is the Dirac measure at $x_0$, integration of $v$
w.r.t. $\mu$ amounts to time integration along the system trajectory starting at $x_0$:%
\[ \langle v,\mu \rangle := \intT \int_{\mathcal{X}} v(t,x)\mu(dt,dx) = \intT v\bigl(t,x(t|x_0)\bigr)\,dt.  \] %

The occupation measure $\mu$ can be disintegrated as $\mu(\mathcal{A}\times\mathcal{B}) = \int_{\mathcal{A}}
\mu_{k}(\mathcal{B}|t)\,dt$ where conditional $\mu_{k}$ is a stochastic kernel, in the sense that for every fixed $t
\in \mathcal{T}$, $\mu_{k}(dx|t)$ is a probability measure on $\mathcal X$ describing the distribution of the state $x$ at
time $t_k$, and for every $\mathcal{B}\in\mathcal{X}$, $t \mapsto \mu_{k}(\mathcal{B}|t)$ is a measurable function on
$\mathcal{T}$.

On the one hand, the introduced measures allow us to consider the whole set of initial conditions ($\mu_{0}$). On
the other hand, it allows a reformulation of the nonlinear dynamics \eqref{eq:system:ode} as a linear equation which only
depends on the occupation measures ($\mu$, $\mu_{t_k}$), as shown next.

\subsection{Liouville's equation}

With these notations, for all sufficiently regular test functions $v \in C^1(\mathcal{T}\times\mathcal{X}; {\mathbb R})$
and $\mathcal T:=[0,1]$, it holds that
\begin{multline}
    \label{eq:bla}
    \int_{\mathcal{X}} v(1,x)\mu_1(dx) - \int_{\mathcal{X}} v(0,x)\mu_0(dx) = 
    \intT \int_{\mathcal{X}} \frac{d}{dt}\,v\bigl(t,x(t|x_0)\bigr)\mu_0(dx_0).
\end{multline}
The right-hand-side of the above equation can be rewritten
\begin{align*}
    &\intT \int_{\mathcal{X}} \Bigl(\frac{\partial}{\partial t}v\bigl(t,x(t|x_0)\bigr) 
+ \mathrm{grad}\:v\bigl(t,x(t|x_0)\bigl) \cdot f\bigl(t,x(t|x_0)\bigr) \Bigr) \mu_0(dx_0)\,dt \\
    =& \intT \int_{\mathcal{X}} \left(\frac{\partial}{\partial t}v(t,x) + \mathrm{grad}\:v(t,x) \cdot f(t,x) \right)
    \mu(dt,dx).
\end{align*}
To simplify notation, we introduce the Liouville operator ${\mathcal L} : {C}^1(\mathcal{T}\times\mathcal{X}) \to
{C}(\mathcal{T}\times \mathcal{X})$ as $ {\mathcal L} v := \frac{\partial v}{\partial t} + \mathrm{grad}\:v \cdot f$ and
its adjoint ${\mathcal L}' : {C}(\mathcal{T}\times\mathcal{X})' \to {C}^1(\mathcal{T}\times\mathcal{X})'$ such that
$\langle {\mathcal L} v,\mu \rangle = \langle v,{\mathcal L}'\mu \rangle$ for all $v \in
{C}^1(\mathcal{T}\times\mathcal{X})$, i.e.  $ {\mathcal L}' \mu := -\frac{\partial \mu}{\partial t} - \mathrm{div}(\mu
f).  $

With these notations, equation~\eqref{eq:bla} can be written concisely as
\begin{align}
    \label{eq:Liouv}
    \langle {\mathcal L} v, \mu \rangle = \langle v,\delta_1 \mu_1 \rangle - \langle v,\delta_0 \mu_0 \rangle
\end{align}
for all $v \in {C}^1(\mathcal{T}\times\mathcal{X})$, where $\delta_0$ and $\delta_1$ refers to $t=0$ and $t=1$,
respectively. Equivalently, in the sense of distributions, we can write
\begin{equation}\label{liouville} {\mathcal L}' \mu = \delta_1\mu_1 - \delta_0\mu_0.
\end{equation}
Equation~\eqref{liouville} is Liouville's equation and is also called the continuity equation in statistical physics or
fluid mechanics.

Whereas the function $x \in {C}(\mathcal{T}\times\mathcal{X})$ satisfies the nonlinear ordinary differential equation
(\ref{eq:system:ode}), the occupation measure $\mu \in {P}(\mathcal{T}\times\mathcal{X})$ satisfies Liouville's equation
(\ref{liouville}), which is a linear partial differential equation (PDE) in the space of probability measures.

Note that as the initial conditions are not known, the measures are unknown
as well.  In the next sections we
derive an optimization problem that allows us to determine the unknown measures.

\subsection{Estimating the region of attraction}\label{s:roa}

In \cite{Henrion_Korda_2013_ROA_occupation_measure} it was proved that the region of attraction $\mathcal{X}_0^\ast$,
defined as the set of initial conditions $x_0$ consistent with the dynamics \eqref{eq:system:ode} and the constraints
$x(t) \in \mathcal{X}$, $t \in \mathcal{T}$ and a constraint $\mathcal{X}_1$ at $t=1$, is the support of the measure
$\mu_0$ solving the infinite-dimensional linear programming (LP) problem
\begin{equation}\label{plp}
    \begin{array}{ll}
        \sup & \langle 1,\mu_0 \rangle \\
        \mathrm{s.t.} & \hat{\mu}_0 + \mu_0 = \lambda,  \\
        & {\mathcal L}' \mu + \delta_0 \mu_0 - \delta_1 \mu_1 = 0, \\
        & \hat{\mu}_0 \geq 0, \: \mu_0 \geq 0, \: \mu_1 \geq 0, \: \mu \geq
        0,
    \end{array}
\end{equation}
where $\lambda$ is the Lebesgue measure restricted to $\mathcal X$, i.\,e.\ the standard $n_x$-dimensional volume.  The
supremum in \eqref{plp} is w.r.t.\ measures $\hat{\mu}_0 \in {P}(\mathcal{X})$, $\mu_0 \in {P}(\mathcal{X})$, $\mu_1
\in{P}(\mathcal{X})$ and $\mu \in {P}(\mathcal{T}\times\mathcal{X})$.  Note that the slack measure $\hat{\mu}_0$ results
from the inequality $\mu_0 \le \lambda$ as further explained in \cite{Henrion_Korda_2013_ROA_occupation_measure}.
The above LP problem has a dual LP
\begin{align}\label{dlp}
        \inf & \langle v_0,\lambda \rangle\nonumber \\
		\mathrm{s.t.} &\begin{array}[t]{lr}
				v_0(x) \geq 0, & \forall x \in \mathcal{X}, \\
        v_0(x) \geq 1+v(0,x), & \forall x \in \mathcal{X}, \\
        v(1,x) \geq 0, & \forall x \in \mathcal{X}_1, \\
        - {\mathcal L} v(t,x) \geq 0, & \forall(t,x)\in \mathcal{T}\times
        \mathcal{X},
    \end{array}
\end{align}
where the infimum is w.r.t. continuous functions $v_0 \in {C}(\mathcal{X})$ and $v \in
{C}(\mathcal{T}\times\mathcal{X})$.

The above  LPs \eqref{plp} and \eqref{dlp} are infinite-dimensional, because the equations are
required to hold for all test functions $v$. One can solve these LPs by a converging hierarchy of finite-dimensional
linear matrix inequality (LMI) problems using semidefinite programming.  At a given relaxation order $d$, the primal LMI
is a moment relaxation of primal LP (\ref{plp}), whereas the dual LMI is a polynomial sum-of-squares (SOS) restriction
of dual LP (\ref{dlp}).

\subsection{Sum-of-squares relaxation of the infinite dimensional LP}
\label{sec:sos}

The dual LMI w.r.t.\ \eqref{dlp} is given by
\begin{equation}
    \label{eq:outer-approx:dual}
    \hspace{-0.5em}%
    \begin{array}{rllll}
        \inf & {v_0^{c}}'l\\
        \text{s.t.}  & -\mathcal{L}v(t,x) &\hspace{-0.75em}=\hspace{-0.75em}& p(t,x) + q_0(t,x)t(1-t)
        +\sum_{i=1}^{m_{x}}q_i(t,x)g_{x,i}(x), \\
        & \hfill v_0(x) &\hspace{-0.75em}=\hspace{-0.75em}&  v(0,x) + 1 + r_0(x) + \sum_{i=1}^{m_x}r_{0,i}(x)
        g_{x,i}(x), \\
        & \hfill v_0(x) &\hspace{-0.75em}=\hspace{-0.75em}& p_0(x) + \sum_{i=1}^{m_{x}}q_{0,i}(x) g_{x,i}(x), \\
        & \hfill v(1,x) &\hspace{-0.75em}=\hspace{-0.75em}& p_1(x) + \sum_{i=1}^{m_{x_k}}q_{1,i}(x) g_{x_1,i}(x),
    \end{array}
\end{equation}
where $l$ is the vector of Lebesgue moments over $\mathcal{X}$ indexed in the same basis in which the polynomial
$v_0(x)$ with coefficients $v_0^c$ is expressed. The minimum is over polynomials $v(t,x)$ and $v_0(x)$, and polynomial
sum-of-squares $p(t,x)$, $q_0(t,x)$, $q_i(t,x)$, $p_0(x)$, $q_{0,i}(x)$, $p_1(x)$, $r_0(x), r_{0,i}(x), \forall i =
\iset{1}{m_x}$ and $q_{1,i}(x), \forall i = \iset{1}{m_{x_k}}$ of appropriate degrees.  The constraints that polynomials
are sum-of-squares can be written explicitly as LMI constraints, and the objective is linear in the coefficients of the
polynomial $v_0(x)$.  Therefore, problem \eqref{eq:outer-approx:dual} can be formulated as an SDP.

From the solution of the dual LMI of order $d$, we obtain a polynomial $v^d_0(x)$ of given degree $2d$ which is such
that the semialgebraic set $\mathcal{O}^d:=\{x_0 \: :\: v^d_0(x) \geq 1\}$ is a valid outer-approximation of the region
of attraction $\mathcal{X}_0^\ast$, i.e.  $\mathcal{X}_0^\ast \subset \mathcal{O}^d$.  Moreover, the approximation
converges in the Lebesgue measure, or equivalently almost uniformly, in the sense that $\lim_{d\to\infty}
\lambda(\mathcal{O}^d) = \lambda(\mathcal{X}_0^\ast)$, see \cite{Henrion_Korda_2013_ROA_occupation_measure} for details.

\section{Consistent Parameter Estimation}
\label{sec:appli}

We use now the occupation measure approach to address the set-based parameter estimation
Problems~\ref{prob:param_estim} and \ref{prob:inconsistency}. As shown next, this requires several extensions to be able to
consider the unknown-but-bounded state constraints at the different measurement time-points $t_k$.

First, we split the solution $\{x(t), \:t \in [0,1]\}$ of problem (\ref{eq:system}) into $m_t-1$ arcs $\{x(t), \:t \in
\mathcal{T}_k\}$ with $\mathcal{T}_k:=[t_k,t_{k+1}]$, $k=0,1,\ldots,m_t-2$. Now consider their respective occupation
measures $\mu_{k,k+1} \in {M}(\mathcal{T}_k\times\mathcal{X})$, with intermediate measures $\mu_{k} \in
{P}(\mathcal{X}_k)$. Obviously $\sum_{k=0}^{m_t-2} \mu_{k,{k+1}} = \mu$ and considering Liouville's equation
(\ref{liouville}) on each arc of the trajectory yields a system of linear PDEs
\[ {\mathcal L}' \mu_{k,k+1} = \delta_{t_{k+1}}\mu_{k+1} - \delta_{t_k}\mu_{k}, \quad k=0,1,\ldots,m_t-2.
\]

\subsection{Outer-approximation}
\label{sec:outer}

An outer-approximation $\mathcal{O}\supseteq \mathcal{X}_0^\ast$ is given by the support of the measure $\mu_0$ solving
the LP
\begin{equation}\label{plpn}
\hspace{-0.6em}
    \begin{array}{ll}
        \sup & \langle 1,\mu_0 \rangle \\
        \mathrm{s.t.} & \hat{\mu}_0 + \mu_0 = \lambda,  \\
        & {\mathcal L}' \mu_{k,k+1} = \delta_{t_{k+1}} \mu_{k+1} - \delta_{t_k}
	\mu_{k}, \\
        & \hat{\mu}_0 \geq 0, \: \mu_0 \geq 0, \\
        &  \mu_{k+1} \geq 0, \: \mu_{k,k+1} \geq 0, \ k=0,1,\ldots,m_t-2, \\
    \end{array}
\end{equation}
where the supremum is w.r.t. measures $\hat{\mu}_0 \in {P}(\mathcal{X})$, $\mu_0 \in {P}(\mathcal{X})$, $\mu_{k+1} \in
{P}(\mathcal{X}_{k+1})$, $\mu_{k,k+1} \in {M}([0,1]\times\mathcal{X})$, $k=0,1,\ldots,m_t-2$.  The above LP problem has
a dual LP
\begin{align}\label{dlpn}
    \inf\ & \langle v_0,\lambda \rangle\nonumber \\
		\mathrm{s.t.} & \begin{array}[t]{lr}
         v_0(x) \geq 0, & \forall x \in \mathcal{X}, \\
         v_0(x) \geq 1+v_{0,1}(0,x), & \forall x \in \mathcal{X}, \\ 
         v_{k-1,k}(t_k,x) \geq v_{k,k+1}(t_k,x), & \forall x \in
	\mathcal{X}_k, k=1,\ldots,m_t-2, \\
         v_{m_t-2,m_t-1}(t_{m_t-1},x) \geq 0, & \forall x \in
						\mathcal{X}_{m_t-1}, \\
         - {\mathcal L} v_{k,k+1}(t,x) \geq 0, &\!\!\!\!\!\!\!\! \forall(t,x)\in
						\mathcal{T}\times\mathcal{X}, k=0,1,\ldots,m_t-2,
    \end{array}
\end{align}
where the infimum is w.r.t. continuous functions $v_0 \in {C}(\mathcal{X})$, $v_{k,k+1} \in {C}([0,1],\mathcal{X})$,
$k=0,1,\ldots,m_t-2$.

As in Section \ref{s:roa}, the above infinite-dimensional LPs (\ref{plpn}) and (\ref{dlpn}) are solved by a converging
hierarchy of finite-dimensional LMI problems. From the solution of the dual LMI corresponding to \eqref{dlpn} we
obtain again a polynomial $v_0^d(x)$ such that $\mathcal O^d:=\{x_0 : v_0^d \geq 1\} \supset \mathcal X^*_0$ and
$\lim_{d\to\infty} \lambda(\mathcal{X}^d_0) = \lambda(\mathcal{X}^\ast_0)$.

\subsection{Inner-approximation}
\label{sec:inner}

For the inner-approximation $\mathcal{I}\subseteq\mathcal{X}_0^\ast$, we build on an idea suggested in
\cite{Korda_etAl_2013__Inner_approximations_occupation_measure}. However, we have to take care of the different
measurements at time-points $t_k$.

In the following, we consider the set of initial conditions $\mathcal{C}^\ast_0$ for which there exists an admissible
trajectory~\eqref{eq:traj_state} that violates at least one of the constraints
that define $\mathcal{X}$ and $\mathcal{X}_k$.
By continuity of solutions (the vector field $f$ is polynomial and Lipschitz on the compact set $\mathcal{X}$), this set
is equal to
\begin{align}
    \label{eq:consistent:X_0:complement:traj}
    \mathcal{C}^\ast_0 \coloneqq \left\{ \hspace{-0.5em}
        \begin{array}{llll}
            x_0  : &\exists x(t|x_0) &\text{s.\,t.\ }\\
            &\ \ \exists t \in [0,1] \text{\ and } \exists \nu &\text{s.\,t.\ }  g_{x,\nu}(x) < 0 \\
            &\text{or\ }\\
            &\ \ \exists t_\kappa \phantom{\in [0,1]} \text{\ and } \exists \eta &\text{s.\,t.\ } g_{x_\kappa,\eta}(x) < 0\ 
        \end{array}
        \hspace{-0.5em} \right\},
\end{align}
where $\nu = \iset{1}{m_{x}}$ and $\eta = \iset{1}{m_{x_k}}$ index the violated constraint, and $\kappa =
\iset{0}{m_t-1}$ describes the time-points. Obviously $\mathcal{X}^\ast_0 \coloneqq \mathcal{X} \setminus
\mathcal{C}^\ast_0$.

Note that depending on $x(t|x_0)$, there are different combinations of constraints ($g_{x_\kappa,\eta}$ and $g_{x,\nu}$)
that can be violated.  We directly deal with the different combinations  and derive for each one
 an outer-approximation of the set of initial conditions (and hence parameters) that lead to the violation
of the constraint.  As will be detailed in the sequel, the inner-approximation is then obtained from the union of complements
of these outer-approximations.

To simplify the presentation, we assume that the constraints defining $\mathcal{X}$ are
not violated, i.\,e.\ $x(t|x_0) \in \mathcal{X}, \forall x_0 \in \mathcal{X}_0, \forall t \in \mathcal{T}$. This is a
mild assumption, since the bounds $\mathcal{X}$ can often be derived from system insight (e.\,g.\ mass conservation in
chemical reaction networks), or can be chosen sufficiently conservative. In any case, the constraints defining
$\mathcal{X}$ can be treated similarly to $\mathcal{X}_k$.

Note that the number of possible combinations, $m_c = 2^{m_t\cdot m_{x_k}}-1$, can be very large.
We can reduce the number of combinations significantly due to the observation that if one constraint at $t_\kappa$ is
violated, then it does not matter if the constraints for $k > \kappa$ are satisfied or not and can therefore be ignored.
This can be formalized by:
\begin{align}
    \label{eq:complements_ij}
    &\hspace{-0.3em}\mathcal{C}_{\kappa,\eta} \coloneqq \\
    &\begin{cases}
        \mathcal{X}_k = \mathcal{X}_k, &\forall k < \kappa\\
        \mathcal{X}_{\kappa} = \left\{
            \begin{array}{rrcl}
                x : &g_{x_\kappa,i}(x) &\ge 0& \forall i \ne \eta,\\
                & g_{x_\kappa,\eta}(x) &< 0
            \end{array}
        \right\}, & \phantom{\forall} k = \kappa\\
        \mathcal{X}_{k} = \mathcal{X}, &\forall k > \kappa
    \end{cases}.\nonumber
\end{align}
where $i = \iset{1}{m_{x_k}}$.
\begin{remark}(Strict and non-strict inequalities) In equation~(\ref{eq:complements_ij}), we consider strict
    inequalities. To account for strict inequalities small
    numbers (slack variables) are introduced when the LMIs are solved.
\end{remark}
Then, once the $m_{x_k}\cdot m_t$ different outer-approximations $\mathcal{O}(\mathcal{C}_{\kappa,\eta})$ have been determined, an inner-approximation is
obtained by
\begin{align}
    \label{eq:inner-approx:compl}
    \mathcal{I} \coloneqq \mathcal{X} \setminus
    \mathop{\bigcup\limits_{\kappa=1,\ldots,m_t-1}}\limits_{\eta=1,\ldots,m_{x_k}}\mathcal{O}(\mathcal{C}_{\kappa,\eta}).
\end{align}

\subsection{Inconsistency certificate}
\label{sec:incon_certificate}

Solving Problem~\ref{prob:inconsistency}, i.\,e.\ certifying emptiness of the set of consistent parameters $\mathcal{X}_0^\ast$, was not addressed
in \cite{Henrion_Korda_2013_ROA_occupation_measure}.
 Mathematically, this amounts to certifying infeasibility of the infinite-dimensional
LP
    \begin{align}
        \mathrm{find} &\ \ \hat{\mu}_0,\mu_k,\mu_{m_t-1},\mu_{k,k+1}\nonumber \\
        \mathrm{s.t.} &\ \ \hat{\mu}_0 + \mu_0 = \lambda,\nonumber  \\
        &\ \ {\mathcal L}' \mu_{k,k+1} = \delta_{t_{k+1}} \mu_{k+1} - \delta_{t_k}
	\mu_{k}, \label{flp}\\
        &\ \ \hat{\mu}_0 \geq 0, \: \mu_0 \geq 0,\nonumber \\
        &\ \ \mu_{k+1} \geq 0, \: \mu_{k,k+1} \geq 0, \quad
	k=0,1,\ldots,m_t-2,\nonumber
    \end{align}
which corresponds to \eqref{plpn} without a cost function.  We can check that we meet all the assumptions to apply the
generalized Farkas theorem of \cite[Theorem 2]{craven} and that non-existence of measures $\hat{\mu}_0,\mu_k,\mu_{m_t-1},\mu_{k,k+1}$
solving LP problem (\ref{flp}) is equivalent to the existence of continuous functions $v_0, v$ solving the dual LP
\begin{align}\label{fdlp}
\mathrm{find} &\ \ \hat{\mu}_0,\mu_k,\mu_{m_t-1},\mu_{k,k+1},
							k=1,\ldots,m_t-2,\nonumber\\
							\mathrm{s.t.} &\begin{array}[t]{lr}
	\langle v_0,\lambda \rangle = -1,\\
			    v_0(x) \geq 0,  &\forall x \in \mathcal{X},\\
					v_0(x) \geq 1+v_{0,1}(0,x), &  \forall x \in
							\mathcal{X},\\
							v_{k-1,k}(t_k,x) \geq v_{k,k+1}(t_k,x),\hspace{-1cm}& \forall x \in \mathcal{X}_k, k=1,\ldots,m_t-2,\\
				v_{m_t-2,m_t-1}(t_{m_t-1},x) \geq 0,  &  \forall x \in
			\mathcal{X}_{m_t-1}, \\
			- {\mathcal L} v_{k,k+1}(t,x) \geq 0,\hspace{-1cm}
			&\!\!\!\!\! \forall(t,x)\in
	\mathcal{T}\times\mathcal{X}, k=0,1,\ldots,m_t-2,
\end{array}
\end{align}

If $\mathcal{X}^\ast_0$ is empty, then LP problem (\ref{fdlp}) is
infeasible. If an LMI relaxation of problem (\ref{fdlp})
is infeasible at some order $d$, certified by a normalized Farkas vector in the dual LMI
relaxation (cf.\ Section~\ref{sec:sos}), LP is also infeasible. Thus, finding a normalized Farkas vector thus gives a sufficient condition that can be used to provide a certificate of inconsistency.


\section{Example}
\label{sec:example}

\begin{figure*}[Ht]
\begin{center}
\includegraphics[width=0.8\linewidth]{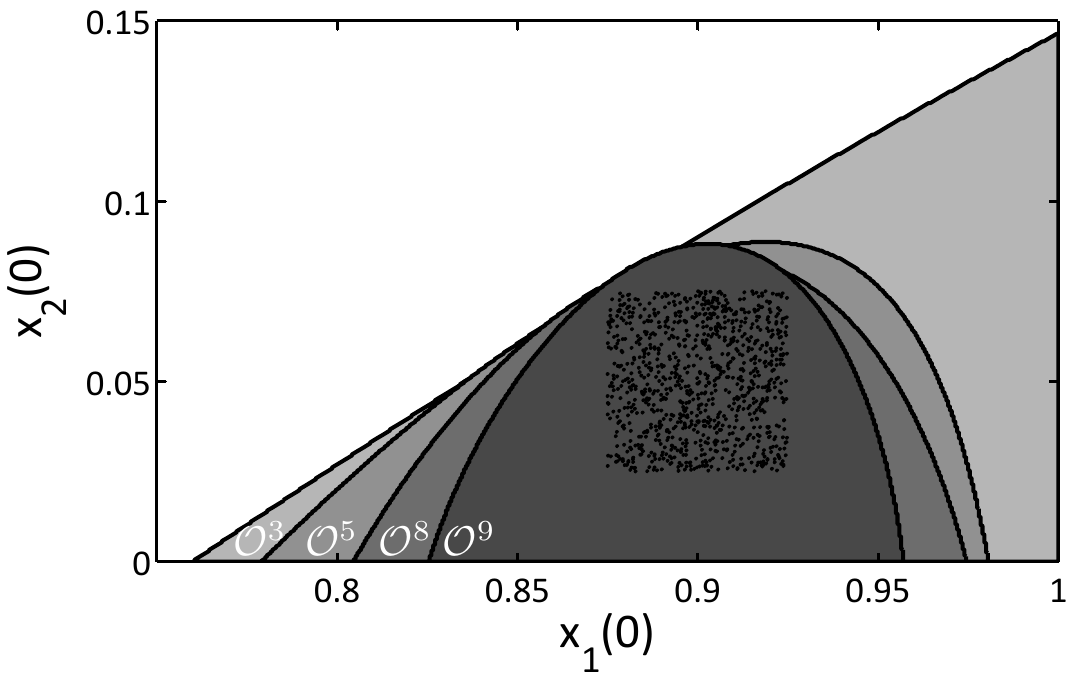}\\
\caption{Set-based estimation of consistent parameters and initial conditions.
 Consistent samples (black dots) were determined by
uniform sampling from $\mathcal{X}$ and subsequent
       numerical integration.
Converging hierarchy of outer-approximations $\mathcal{O}^d$ with $d=3,5,8,9$
for estimation of $x_1(0)$, $x_2(0)$. \label{fig:plotsa}}
\end{center}
\end{figure*}

\begin{figure*}[Ht]
\begin{center}
\includegraphics[width=0.8\linewidth]{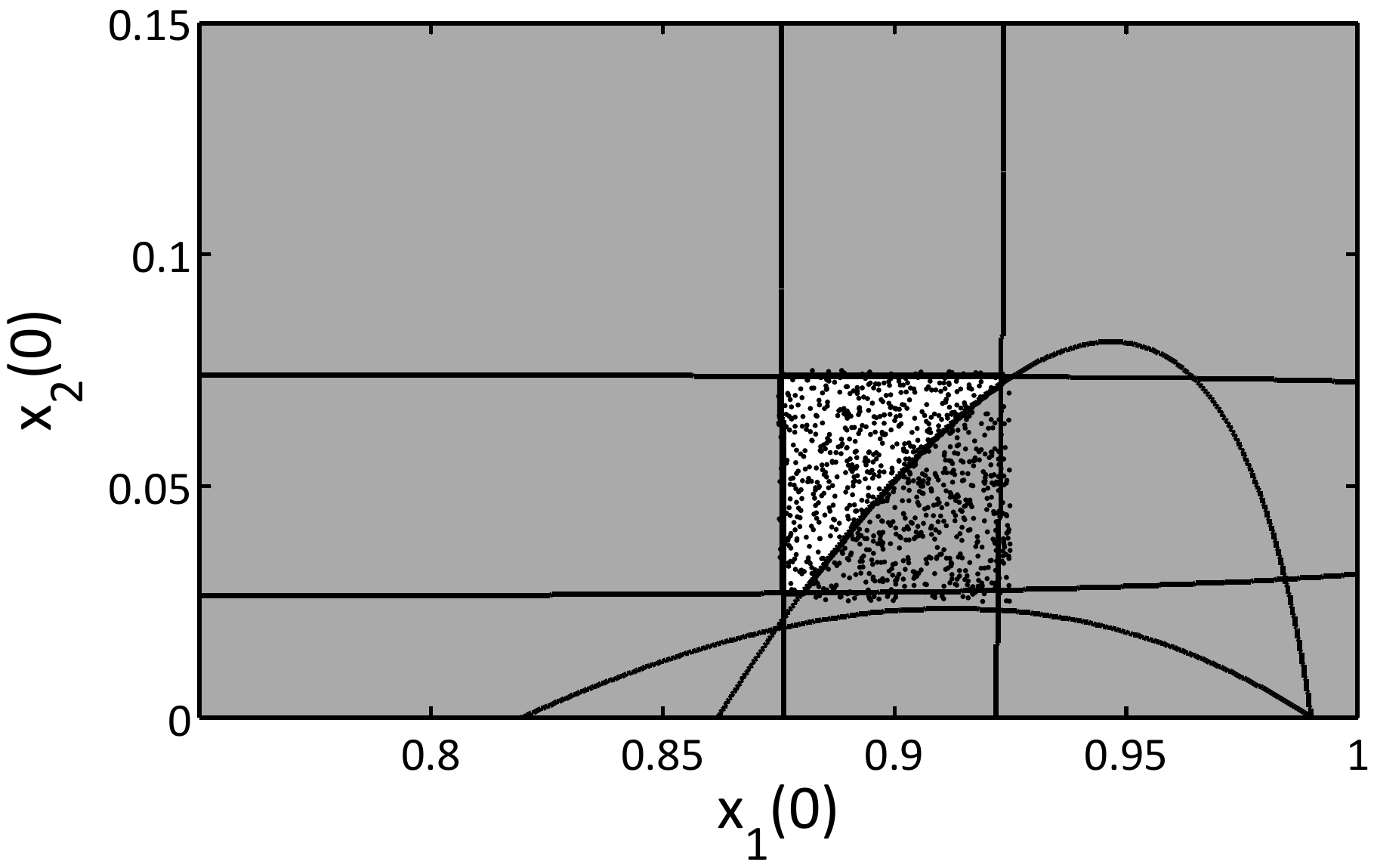}\\
\caption{Set-based estimation of consistent parameters and initial conditions.
 Consistent samples (black dots) were determined by
uniform sampling from $\mathcal{X}$ and subsequent
       numerical integration.
Inner-approximations (white area) results from the union of outer-approximations
(grey area and black lines) for $d=7$. \label{fig:plotsb}}
\end{center}
\end{figure*}

\begin{figure*}[Ht]
\begin{center}
\includegraphics[width=0.8\linewidth]{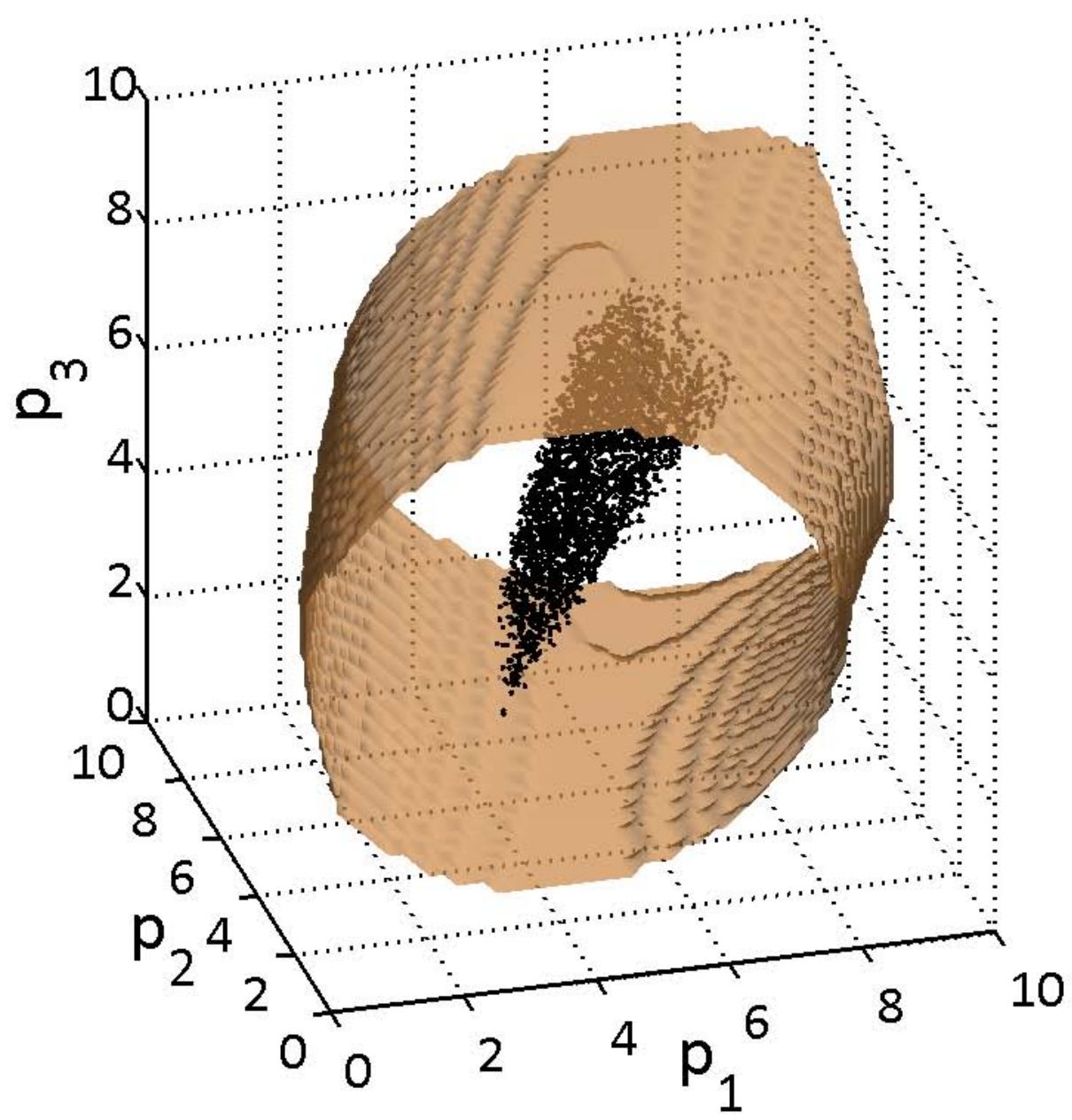}\\
\caption{Set-based estimation of consistent parameters and initial conditions.
 Consistent samples (black dots) were determined by
uniform sampling from $\mathcal{X}$ and subsequent
       numerical integration.
Outer-approximations ($d=3$) of consistent initial conditions $x_1(0)$, $x_2(0)$ and parameters $p_1, p_2, p_3$ projected numerically onto the parameter space. \label{fig:plotsc}}
\end{center}
\end{figure*}

We consider a biochemical reaction network in which the substrate $x_1$ is enzymatically converted into a product via
intermediary complex $x_2$ \cite{schnell2006mechanism}. The continuous-time dynamics is given by:
\begin{equation}
    \label{eq:example}
\hspace{-.5cm}
    \begin{array}{rcl}
        \dot{x}_1(t) &=& -p_1\, x_1(t) \bigl(1-x_2(t)\bigr)  + p_2\,
			x_2(t), \\
        \dot{x}_2(t) &=& \phantom{-}p_1\, x_1(t)\bigl(1-x_2(t)\bigr)
			-\bigl(p_2+p_3\bigr)\,x_2(t).
    \end{array}
\end{equation}

The following constraints on the parameters and states were used:
\begin{equation*}
    \mathcal{X} : \left\{
        \begin{array}{rcll}
            (p_i-0.1)(10-p_i) &\geq& 0, &i = 1,2,3, \\
            x_i(1-x_i) &\geq& 0, &i=1,2.
        \end{array}
    \right.
\end{equation*}
Measurement constraints $\mathcal{X}_k$ were generated from a simulated nominal trajectory with $p_1=p_2=p_3=5.05,
x_1(0) = 0.9$ and $x_2(0) = 0.05$. An absolute uncertainty of $\pm 0.025$ was added to the nominal values at the
measurement time points $t_0 = 0, t_1 = 0.3, t_2 = 1$ to obtain $\mathcal{X}_0$, $\mathcal{X}_1$, $\mathcal{X}_2$. Thus,
at each time-point four constraints were used, which results in $m_c = 4\cdot 3$.  To avoid numerical troubles, note
that the dynamics \eqref{eq:example} and $\mathcal{X}$ were scaled such that the values of the parameters range from 0 to 1.

As can be seen in Figure~\ref{fig:plotsa}, the sequence of outer-approximations $\mathcal{O}^d$ converges to the
consistent parameter set for increasing $d$. Using YALMIP and Sedumi, the solving time was about 5 seconds for $d=3$, and 40
minutes for $d=9$.

The inner-approximation in Figure~\ref{fig:plotsb} corresponds to the complement (in $\mathcal{X}$) of the union of
twelve outer-approximations $\mathcal{O}(\mathcal{C}_i)$. However, only eight outer-approximations are shown since the
other four were empty.  Solving time was on average 5 minutes per problem.

In addition we determined the outer-approximation of consistent initial conditions $x_1(0), x_2(0)$ and parameters $p_1,
p_2, p_3$, see Figure~\ref{fig:plotsc}. Solving time was about 4 hours.

Note that the inconsistency certificates derived in Section~\ref{sec:incon_certificate} could also be used to invalidate
entire regions in the space of the parameters and initial conditions.

\section{Discussion and Conclusions}
\label{sec:conclusion}

Occupation measures are a classical tool in kinetic theory, statistical physics, optimal transport, Markov decision
processes, amongst others. In a broad perspective, the potential of occupation measures and subsequent convex
relaxations are not often used in systems control.  We used occupation measures to approximate consistent parameter sets
for continuous-time nonlinear systems without the need of numerical integration. A particular feature of the derived
approximations is the (almost uniform) convergence to the true (possibly nonconvex) consistent parameter set.  As
demonstrated at the example, outer- and inner-approximations can  be obtained even though only few measurements with
relatively large error are used. Tighter approximations are expected if more measurements are used, or if e.\,g.\
outer-approximations are iteratively used to refine the results as it proved useful for linear relaxations \cite{Streif_etAl_2012_ADMIT,Rumschinski_etAl_2010_Set-based_parameter_estimation_model_invalidation,Borchers_etAl_2009_CDC__Set-based_estimation}.

While inconsistency certificates can be used to prove inconsistency of entire models or parameter regions, the outer-
and inner-approximations can be used to get the shape of the consistent parameter set.  Such a description of the
outer-approximation by a single polynomial can be very useful in some applications. In other cases different types of
representations like a collection of half-spaces might be more beneficial.  The inconsistency certificates could also be
used to derive in an iterative and recursive manner either outer bounding boxes or a description of the consistent set
using a bisection algorithms (cf.\ \cite{Streif_etAl_2012_ADMIT}).

A drawback of the presented approach is the computational workload resulting from the LMI constraints. Theoretically,
the resulting problems can be solved in polynomial time w.r.t.\ the input size.  However, due to the size limitations of
state-of-the-art SDP solvers, this approach is at the moment restricted to problems of small dimensions.
An alternative could be the use of LP relaxations for larger dimensional systems. Note that the specific geometry of the
LMI constraints makes these relaxations typically much more accurate than the LP relaxations and a trade-off between
accuracy and problem size has to be made, see in particular the discussion in \cite[Section 5.4.2]{Lasserre_2009}.

As an interesting extension, one could consider statistical information, i.\,e.\ further constraints on the moments of
the occupation measures. Here we assumed no more information such as statistics or probability distributions to be
given. In many applications where there is a limited number of replicates, i.\,e. too few to obtain a meaningful
statistic, this is actually the case. However, statistical information can be included if the data are polynomial or
information about moments are available
\cite{Lasserre_2009,Savorgnan_etAl_2009_CDC__Optimal_Control_occupation_measures_moment_relaxations}.


\end{document}